\documentclass{article}

\title{On a $p$-adic extension of the Jacquet-Langlands correspondence to weight~1}
\author{L. J. P. Kilford}

\usepackage{amsmath,amssymb,ifthen,amsthm,url}

\newcommand{\abcds}{\left(\begin{smallmatrix}a & b \\ c & d\end{smallmatrix}\right)}
\newcommand{\abcdsi}{\left(\begin{smallmatrix}a_i & b_i \\ c_i & d_i\end{smallmatrix}\right)}
\newcommand{\smatrix}[4]{\left(\begin{smallmatrix}#1 & #2 \\ #3 & #4\end{smallmatrix}\right)}

\newcommand{\Pari}{{\sc Pari}}

\newcommand{\Magma}{{\sc Magma}}

\newcommand{\Q}{\mathbf{Q}}
\newcommand{\Z}{\mathbf{Z}}
\newcommand{\A}{\mathbf{A}}
\newcommand{\C}{\mathbf{C}}
\newcommand{\N}{\mathbf{N}}
\newcommand{\M}{\mathcal{M}}
\newcommand{\m}{\mathfrak{m}}
\newcommand{\T}{\mathbf{T}}

\newtheorem{theorem}{Theorem}

\newtheorem{definition}[theorem]{Definition}

\DeclareMathOperator{\GL}{GL}

\DeclareMathOperator{\SL}{SL}

\newcommand{\Gouvea}{Gouv\^ea}

\begin{document}
\maketitle

\begin{abstract}
In this paper, we consider a novel version of the classical Jacquet-Langlands \mbox{correspondence}, explore a $p$-adic extension of the correspondence, and as an explicit example we find an overconvergent automorphic form of weight~1 which corresponds to a classical modular form of weight~1, using both experimental and theoretical methods.
\end{abstract}

\section{Introduction}

The Langlands Program has been one of the major organizing programs of the $20^{\rm th}$ century. It seeks to relate number theoretic and representation theoretic objects by attaching L-functions to both sides and then showing that these are equal. The first, motivating, example was the Artin reciprocity law, which generalizes quadratic reciprocity. This was followed by work of Hecke which associated Dirichlet L-functions to automorphic representations. An introductory reference for the Langlands Program is~\cite{bump-et-al}.

This paper uses the theory of overconvergence, both for automorphic forms and for modular forms. This has been used before to prove results about classical modular forms; for instance, in the work of Buzzard-Kilford~\cite{buzzard-kilford}, Emerton~\cite{emerton-thesis}, Jacobs~\cite{jacobs-thesis}, the author~\cite{kilford-2slopes} and Smithline~\cite{smithline-published}, the theory of overconvergence is used to shed new light on classical modular forms. Some interesting computations with $p$-adic modular forms were also performed in~\cite{coleman-stevens-teitelbaum} and~\cite{gouvea-mazur}.

This leads into another motivation for studying this issue; our knowledge of classical modular forms of weight~1 is incomplete, because one cannot compute them using the standard algorithms. The computation of weight~1 forms in general is still an ongoing subject of research; see for example~\cite{edixhoven-mod-p}; although \Magma{} will now compute modular forms of weight~1, the algorithm used is not as efficient as that for higher weight. Indeed, even finding the dimension of spaces of weight~1 forms is difficult; see~\cite{duke-dimension} for a nontrivial bound on the dimension. On the other hand, once we have set up the framework for overconvergent automorphic forms, weight~1 forms are no harder to study than those of any  other weight, so we can use information by studying automorphic forms to prove results about modular forms.

This paper will also give an explicit algorithm for computing the action of the Hecke operators on spaces of automorphic forms, and present some specific examples where we have done this.

We will first recall the definition of classical automorphic forms that we will use in this paper, to fix our notation.

\subsection{Automorphic Forms}
This section follows Section~4 of~\cite{buzzard-families} in notation and approach; this reference presents the subject in a detailed and clear manner.

Let~$p$ be a prime number, and let~$D$ be a definite quaternion algebra over~$\Q$ with discriminant~$\delta$ prime to~$p$. Let~$\mathcal{O}_D$ be a fixed maximal order of~$D$; we also fix an isomorphism~$\mathcal{O}_D \otimes B \cong M_2(B)$, where~$B:=\lim_{\leftarrow}(\Z/M\Z)$, where we take the limit over all integers~$M$ prime to~$\delta$. This isomorphism induces isomorphisms $\mathcal{O}_D \otimes \Z_l \cong M_2(\Z_l)$ and~$\mathcal{O}_D \otimes \Q_l \cong M_2(\Q_l)$ for all primes~$l \nmid \delta$; we will identify these rings with each other.

Let~$\A_f$ be the finite adeles over~$\Q$; we define~$D_f = D \otimes_\Q \A_f$; this can be thought of as the restricted product over all primes~$l$ of~$D \otimes \Q_l$; if~$g \in D_f$ then the component~$g_p$ of~$g$ at the prime~$p$ can be viewed as an element of~$M_2(\Q_p)$.

Let~$U$ be an open compact subgroup of~$D^\times_f$. It is well known that one can write~$D^\times_f$ as a \emph{finite} union of disjoint double cosets, of the form
\[
D^\times_f = \coprod_{i \in I} D^\times d_i U.
\]
We will see later how to compute these~$d_i$.

If~$M$ is a positive integer, then we define~$U_1(M)$ to be the open compact subgroup of~$D^\times_f$ whose elements~$g$ have component~$g_p=\abcds$ at~$p$ with~$c \equiv 0 \mod p$ and~$d \equiv 1 \mod p$, and can be arbitrary at all other places. If~$G$ is a subgroup of~$D^\times$ of finite index, we also define the open subgroup~$U_1(M) \cdot G$ to be the subgroup of~$U_1(M)$ whose image at primes dividing~$\delta$ is in~$G$.

If~$\alpha \ge 1$, we define~$M_\alpha$ to be the monoid consisting of 2 by 2 matrices $\smatrix{a}{b}{c}{d}$ over~$\Z_p$ with nonzero determinant such that~$p^\alpha | c$ and~$p \nmid d$.

Let~$K$ be a complete subfield of~$\C_p$ (the $p$-adic completion of~$\Q_p$). Let~$L_k$ be the space of polynomials over~$K$ in one variable~$z$ of degree at most~$k-2$. We will equip this with an action of~$M_\alpha$; let~$\gamma:=\smatrix{a}{b}{c}{d}$ and let~$h \in L_k$. Then we define the right action by
\[
(h|_{\gamma})(z):=(cz+d)^{k-2} \cdot h\left(\frac{az+b}{cz+d}\right).
\]

We can now define classical automorphic forms; this definition is given in~\cite{buzzard-families}, page~33.
\begin{definition}
Let~$M$ be a positive integer which is prime to~$\delta$ and let~$G$ be a subgroup of~$D^\times$ of finite index. We define $S^D_k(U_1(M))$ to be the space of classical automorphic forms of level~$U_1(M)$ and weight~$k$ for~$D$, which is 
\[
\left\{ f: D^\times_f \rightarrow L_k:\;f(dgu)=f(g)u_p \text{ for all }d\in D^\times,\;u \in U_1(M) \right\}
\]
Similarly, we define $S^D_k(U_1(M) \cdot G)$ to be
\[
\left\{ f: D^\times_f \rightarrow L_k:\;f(dgu)=f(g)u_p \text{ for all }d\in G,\;u \in U_1(M) \right\}.
\]
\end{definition}
Elements of~$S^D_k(U_1(M))$ or~$S^D_k(U_1(M) \cdot G)$ are sometimes called \emph{quaternionic modular forms} in the literature.

Because we can write~$D^\times_f$ as a finite union, we can determine an automorphic form~$f$ by its values on the~$\{d_i\}$; in other words, by the tuple of polynomials $\{f(d_1),\ldots,f(d_n)\}$. We will use this later to perform calculations.

We can define Hecke operators on~$S^D_k(U_1(M))$ as double coset operators in the following way. Let~$\eta \in D^\times_f$ be an element of~$M_\alpha$, and if we have~$f:D^\times_f \rightarrow A$ then we define the right action $f|_{\eta}:D^\times_f \rightarrow A$ by
\begin{equation}
\label{right-action}
(f|_\eta)(g) = f(g\eta^{-1})\cdot \eta_p.
\end{equation}
We now consider the double coset~$U\eta U$; this can be written as a finite union $\coprod_{i}U \eta_i$, with effectively computable~$\eta_i$. We define the Hecke operator~$[U\eta U]$ by
\begin{eqnarray}
[U\eta U]: S^D_k(U_1(M)) &\rightarrow& S^D_k(U_1(M))\\
{[U\eta U]}(f)&=&\sum_i f|_{\eta_i};
\end{eqnarray}
we will see later how this can be computed and exhibit some results of our computations.

There are certain Hecke operators that we will concentrate on; these are the analogues of~$T_p$ and~$U_p$ in the classical setting. We let~$l$ be a prime which does not divide~$\delta$, and define~$\omega_l \in \A_f$ to be the finite adele which is~$l$ at the place~$l$ and is the identity at all of the other finite places. This can be viewed as an element of~$D^\times_f$ via the diagonal embedding. We define~$\eta_l:=\smatrix{\omega_l}{0}{0}{1}$ to be the element of~$D^\times_f$ which is $\smatrix{l}{0}{0}{1}$ at the place~$l$ and the identity at all of the other places. We define~$T_l:=[U\eta_l U]$ (it will be clear from context whether we mean the Hecke operator on classical or automorphic forms). We will call the characteristic polynomial of a Hecke operator the Hecke polynomial of that operator.

We now fix notation for classical (elliptic) modular forms. Let~$\Gamma$ be a congruence subgroup of level~$N$ and let~$\chi$ be a Dirichlet character. We define~$S_k(\Gamma,\chi)$ to be the vector space of classical modular forms of weight~$k$, level~$\Gamma$ and \mbox{character}~$\chi$. This is equipped with the standard Hecke operators~$T_q$ and~$U_p$. There are many standard books which give an introduction to the theory of classical modular forms; for instance, see~\cite{diamond-shurman}.

We now cite a standard version of the correspondence between classical automorphic forms, as we have defined them above, and elliptic modular forms. This can be derived from Theorem~16.1 of~\cite{jacquet-langlands}.

\begin{theorem}[Jacquet-Langlands, Shimizu, Arthur]
\label{standard-jacquet-langlands}
Let~$k \ge 3$ be an integer and let~$M$ be a positive integer prime to~$\delta$. There is an isomorphism between the spaces~$S^D_k(U_1(M))$ and $S^{\gamma-new}_k(\Gamma_1(M)\cap\Gamma_0(\gamma))$ which commutes with the action of the Hecke operators defined above. If~$k=2$ then~$S^{\gamma-new}_2(\Gamma_1(M)\cap\Gamma_0(\gamma))$ is isomorphic to the quotient of~$S^D_k(U_1(M))$ by the subspace of forms which factor through the norm map, and this isomorphism also commutes with the action of the Hecke operators.
\end{theorem}

\section{Families of modular and automorphic forms}

Let~$p$ be a prime number and let~$k_0$ be a non-negative integer. We define a \emph{($p$-adic) family of modular forms} to be a set~$\{f_i\}_{i \in \N}$ of modular forms, where~$f_i$ has weight~$k_0+(p-1)\cdot p^i$ which satisfy the congruence
\[
f_i(q) \equiv f_0(q) \mod p^i.
\]
An example of a $p$-adic family of modular forms is given by the set of Eisenstein series~$\{E_{k+(p-1)p^r}\}_{r \in \N}$ for any positive integer~$k$; this is called the ``Ur-example'' of a $p$-adic family in~\cite{coleman-mazur}.

The $p$-adic limit of such a family was defined by Serre 
in~\cite{serre-mod-p-modular-forms} to be a 
\emph{$p$-adic modular form}; we see that classical modular forms are automatically modular forms under this definition. A more modern description of $p$-adic modular forms can be found in~\cite{gouvea-thesis}; it can be shown that taking the limit either with respect to the weight or to the level gives the same space of modular forms.

A good example of a $p$-adic modular form which is not a classical modular form is given by~$E_2$, which has Fourier expansion given by
\[
E_2(q):=1-24\sum_{n=1}^\infty \sigma_1(n) q^n.
\]
It is well known that~$E_2$ is not a classical 
modular form --- it 
fails to transform properly under the action of~$\SL_2(\Z)$, because the series that defines it is not absolutely convergent --- but it can be shown that it is a $p$-adic modular form for every prime~$p$, because we can $p$-adically approximate it arbitrarily well with classical eigenforms.

Let~$f_A$ be an automorphic form of weight~$k_0$. We define a \emph{family of automorphic forms} to be a set of automorphic eigenforms~$\{f_{A,i}\}_{i \in \N}$ of weight~$(p-1)\cdot p^i$ whose Hecke eigenvalues are congruent to those of~$f_A$ modulo~$(p-1)\cdot p^i$.

One problem with the space of $p$-adic modular forms of a given weight is that~$U_p$ is not a compact operator on this space. We will now introduce a subspace of the $p$-adic modular forms on which~$U_p$ is compact, the overconvergent modular forms.

\section{Overconvergent forms}
In this section, we will give definitions for the rings of overconvergent modular forms and overconvergent automorphic forms.

We first give the definition of overconvergent modular forms. For more details of the construction, see Section~2 of~\cite{coleman-mazur}.
\begin{definition}[Coleman,~\cite{coleman-overconvergent}, page~397]
Let~$p$ be a prime number and let~$m$ be a positive integer. Let~$w$ be a rational number such that~$0 < w < p^{2-m}/(p+1)$. 

We think of~$X_0(p^m)$ as a rigid space over~$\mathbf{Q}_p$, and we let~$t \in X_0(p^m)(\overline{\mathbf{Q}}_p)$ be a point, corresponding either to an elliptic curve defined over a finite extension of~$\Q_p$, or to a cusp.

We define~$Z_0(p^m)(w)$ to be the connected component of the affinoid
\[
\left\{t \in X_0(p^m): \; v_p(E_4(t)) \le w\right\}
\]
which contains the cusp~$\infty$.

Let~$\mathcal{O}$ be the structure sheaf of~$Z_0(p^m)(w)$.
We call sections of~$\mathcal{O}$ on the affinoid~$Z_0(p^m)(w)$
\emph{$w$-overconvergent $p$-adic modular forms of weight~$0$ and level~$\Gamma_0(p^m)$}.
If a section~$f$ of~$\mathcal{O}$ is a $w$-overconvergent modular form, then we say that~$f$ is an \emph{overconvergent $p$-adic modular form}. We denote the space of $w$-overconvergent $p$-adic modular forms of weight~0 by~$\M_0(p^m,w)$.

We now let~$\chi$ be a primitive Dirichlet character of conductor~$p^m$ and let~$k$ be an integer such that~$\chi(-1)=(-1)^k$. Let~$E^*_{k,\chi}$ be the normalized Eisenstein series of weight~$k$ and character~$\chi$ with nonzero constant term. The space of $w$-overconvergent $p$-adic modular forms of weight~$k$ and character~$\chi$ is given by
\[
\M_{k,\chi}(p^m,w):=E^*_{k,\chi} \cdot \M_0(p^m,w).
\]
\end{definition}
The canonical non-example of an overconvergent modular form is the $p$-adic modular form~$E_2$. In~\cite{koblitz-e2} and~\cite{coleman-gouvea-jochnowitz} it is shown that~$E_2$ is never a $w$-overconvergent modular form, for any choice of~$p$ and~$w$. In fact the latter paper proves that $E_2$ is transcendental over the ring of overconvergent modular forms.

We can show that this is a Banach space over~$K$, and it is also the 
case that classical modular forms are automatically overconvergent 
modular 
forms; we will see later that they can be picked out from the space 
of overconvergent forms. They form a subring of the ring of $p$-adic 
modular forms. 

It is possible to define the space of overconvergent forms of 
weight~$k$ and character~$\chi$ directly, but for applications it is 
convenient to work out the theory for weight~0 forms and then multiply 
by a suitable Eisenstein series of the correct weight and character to 
obtain the weight~$k$ and character~$\chi$ space. One can find the 
action of the~$U_p$ operator in this case also; 
see~\cite{kilford-2slopes} and~\cite{buzzard-kilford} for more details 
in the specific case where~$p=2$.

It should also be noted that, although the definition given here depends on the overconvergence parameter~$w$, the characteristic power series of the~$U_p$ operator is independent of the choice of~$w$. This means that we can effectively perform computations for one choice of~$w$, and know that the results we obtain are independent of the choice of~$w$.

We will now consider overconvergent automorphic forms. We will first generalize the space~$L_k$ to more general, and infinite-dimensional, rings. 

\begin{definition}
Let~$K$ be a complete subfield of~$\C_p$. Then we define~$A_{k,1}$ to be the ring~$K\langle z \rangle$ of power series~$\sum_{n \in \N} a_n z^n$ such that~$a_n \rightarrow \infty$ as~$n \rightarrow \infty$. We call these \emph{convergent power series}.
\end{definition}
We see that the right action by elements of~$M_\alpha$ on~$A_{k,1}$, extending the right action of~$M_\alpha$ on~$L_k$,
\[
(h|_{\gamma})(z):=(cz+d)^{k-2} \cdot h\left(\frac{az+b}{cz+d}\right)
\]
does send~$A_{k,1}$ to itself, so this action is well-defined.

We can now define \emph{overconvergent automorphic forms}.
\begin{definition}
Let~$M$ be a positive integer which is relatively prime to~$\delta$. We define $S^{D,\dagger}_k(U_1(M))$ to be the space of \emph{overconvergent automorphic forms} of level~$U_1(M)$ and weight~$k$ for~$D$, which is 
\[
\left\{ f: D^\times_f \rightarrow A_{k,1}:\;f(dgu)=f(g)u_p \text{ for all }d\in D^\times,\;u \in U_1(M) \right\}.
\]
We can define overconvergent automorphic forms for other open compact subgroups in a similar way.
\end{definition}
We see that these spaces are infinite-dimensional, and also that they contain the spaces of classical automorphic forms, because a polynomial is certainly a convergent power series.  On the other hand, we can determine an overconvergent form~$f$ by its values on the~$\{d_i\}$; in other words, by the tuple of convergent power series~$\{f(d_1),\ldots,f(d_n)\}$.

We also note that the action of the Hecke operators~$T_l$ and~$U_p$ is well-defined; we use the same definitions as we used in the classical case, extending the weight~$k$ action from polynomials to convergent power series.

It is emphasized in the introduction to~\cite{buzzard-families} that one advantage of considering overconvergence for automorphic forms is that the geometric aspects of the definition are much simpler than those in the definition of overconvergent modular forms; we simply change the base ring to define overconvergent automorphic forms, whereas we need to set up a lot of machinery to define overconvergent modular forms.

\section{Extending the classical correspondence}

The Jacquet-Langlands correspondence as it is usually stated gives us a one-to-one correspondence between elliptic modular forms and classical automorphic forms. However, we will be considering an explicit example of a modular form which does \emph{not} lie in the spaces of modular forms considered in Theorem~\ref{standard-jacquet-langlands}.

In this section, we will show how to obtain a version of the classical Jacquet-Langlands correspondence that will be more useful for our purposes. Kevin Buzzard contributed the ideas used in this section.

We will first need to set some notation up. Let $R$ be the ring of Hurwitz integers in $D$, that is the elements of~$D$ whose norms are in $\Z_2$, and let $\m$ denote the maximum two-sided ideal of~$R$, which is the ring of elements of~$D$ with norm  in $2\Z_2$. Then $1+\m \subset R^* \subset D^*$ and $1+m$ can be shown to be normal in $D^*$. 
\begin{theorem}
\label{jacquet-langlands-extended-classical}
Let~$k \ge 3$ be an integer and let~$M$ be an odd positive integer. Let~$D$ be the non-split quaternion algebra over $\Q_2$, which has a norm map $N:D\rightarrow \Q_2$ given by~$N(\alpha)=\alpha\overline{\alpha}$. 

There is an isomorphism 
\begin{equation}
\label{new-jacquet-langlands-relation}
S^D_k(U_1(M)\cdot (1+\m))\cong \left(S^{4-{\rm new}}_k(\Gamma_1(M)\cap\Gamma_0(4))\right)^2 \oplus S^{2-{\rm new}}_k(\Gamma_1(M)\cap\Gamma_0(2)).
\end{equation}
This isomorphism commutes with the action of the Hecke operators defined above. If~$k=2$ then we must quotient out the left hand side of~\eqref{new-jacquet-langlands-relation} by the norm map to obtain an isomorphism. 
\end{theorem}
We note that we can obtain other isomorphisms between modular forms and automorphic forms, if we take still smaller subsets of the unit group of~$O_D$; in this case, we are using the subgroup~$1+\m$ of~${O}^\times_D$. These will correspond to different spaces of modular forms with larger powers of~2 in the level.
\begin{proof}

Let~$\pi$ be a smooth irreducible representation of $D^*$ and let~$v \in \pi$ be a nonzero element which is fixed by~$1+\m$. Let $V$ denote the space of $1+\m$-fixed vectors. We can see immediately that~$V \ne 0$ as $v \in V$, that $V$ is finite-dimensional by smoothness results, and that it is $D^\times$-invariant because~$1+\mathfrak{m}$ is a normal subgroup of~$D^\times$. This means that $V=\pi$ is an irreducible finite-dimensional representation of $D^*/(1+\m)$.

One can check explicitly that this group has two generators, which we can call $a:=(1+i+j+k)/2$ and $b:=1+i$ which satisfy~$a^3=1$ and~$bab^{-1}=a^{-1}$. The eigenvalues of~$a$ on~$V$ will all be cube roots of unity
and if $\zeta$ is a primitive cube root of unity then $\zeta$ and $\zeta^{-1}$ will occur with
equal multiplicity, as~$a$ is conjugate to its inverse. 

Let $V=V_1+V_2$ with
$a=1$ on $V_1$ and $a$ acting by $\zeta$ or $\zeta^{-1}$ on $V_2$. Then $V_1$ and $V_2$ are
both $b$-invariant and hence either $V=V_1$ or $V=V_2$. If $V=V_1$ then any
eigenvector for $b$ is an invariant subspace so $V$ is 1-dimensional. If $V=V_2$
then $V$ is at least 2-dimensional and again if $w$ is an eigenvector for $b$
then ${w,aw}$ span an invariant subspace (as $a^2=-a-1$ on $V_2$) so $V$ is
2-dimensional.

The two cases can be distinguished by whether $a$ acts trivially or not and
this is exactly whether there are $R^*$-fixed vectors or not. Hence we 
see that there is an injection from the
$R^*$-level structure space to the $1+\m$-level structure space but, 
contrary to the classical case, the quotient is all new, and all new 
eigenvalues appear with multiplicity two.

We note that the multiplicities in the right-hand side of~\eqref{new-jacquet-langlands-relation} are the opposite of those in the classical case, where newforms appear with multiplicity one and the oldforms appear with higher multiplicity.
\end{proof}
\section{Computing the Hecke polynomials}
\label{computations}
In the thesis of Jacobs~\cite{jacobs-thesis}, an explicit computation of the~$d_i$ is given. We will repeat this construction here, and generalize it.

For the remainder of this paper we will assume that~$D$ is the Hamiltonian quaternions, which have discriminant~2, and that the prime number~$p$ is greater than~2. The constructions we will give can be generalized to the case when~$\delta=2$.

Let~$U=U_1(p^n)$. This is a compact open subgroup of~$D^\times_f$. Then there is a sequence of isomorphisms of groups given by
\begin{eqnarray}
D^\times\backslash D^\times_f / U &=& D^\times \backslash U_0(1)/U \\
&=& D^\times \cap U_0(1)\backslash U_0(1)/U\\
&=& \mathcal{O}^\times_D\backslash U_0(1)/U\\
&=& \mathcal{O}^\times_D \backslash\SL_2(\Z/p^n\Z)/H,
\end{eqnarray}
where
\[
H:=\{\abcds \in SL_2(\Z/p^n\Z) : \; \abcds \equiv \smatrix{1}{\star}{0}{1} \mod p^n\}.\]
These are all consequences of elementary group theoretic isomorphisms, and the useful fact that~$D^\times_f = D^\times U_0(1)$, where~$U_0(1)$ is the compact subgroup of~$D^\times_f$ which has no restriction at any place~$p$.

It can be shown that, if
\[
G.x:=\left\{\begin{pmatrix} x \\ y \end{pmatrix}:\;x,y \in \Z/p^n\Z \;|\; \text{ at least one of }x,y\in (\Z/p^n\Z)^\times\right\},
\]
then we have
\[
D^\times \backslash D^\times_f / U = \mathcal{O}_D^\times \backslash G.x,
\]
where~$\mathcal{O}_d^\times$ acts on~$G.x$ by matrix multiplication on the left. If~$m$ is the number of orbits then we have the following equality:
\[
\mathcal{O}_D^\times \backslash G.x = \coprod_{i=1}^m \mathcal{O}_D^\times \cdot s_i,
\]
where~$s_i$ is a column vector representing the orbit. Also, if we wish to add extra level structure at~2, we can change the group~$\mathcal{O}_D^\times$ to a subgroup of the units; for instance, we could have~$1+\m$, where~$\m$ is the maximal two-sided ideal of~$\mathcal{O}_D$.

We can compute these~$s_i$ effectively with \Pari{}, and we find the~$d_i$ from them by lifting them to~$\GL_2(\Z_p)$; then~$d_i$ is this lift of~$s_i$ at~$p$, and trivial at the other places.

As an explicit example of this, let~$U=U_1(7)$. We can compute the number of cosets~$m$ to be~2, and suitable~$s_i$ are~$s_1=\left(\begin{smallmatrix} 0 \\ 1\end{smallmatrix}\right)$ and~$s_2=\left(\begin{smallmatrix} 1 \\ 4\end{smallmatrix}\right)$. We can then lift these to entries of~$U_1(7)$ by lifting the~$s_i$ to matrices over~$\GL_2(\Z_7)$ and setting their values at odd primes away from~7 to be trivial.

We now write down the action of the Hecke operators explicitly, in a way that allows us to perform calculations. We let~$U_p:=[U\gamma U]$ be the Hecke operator that we defined earlier; if~$f \in S^D_k(U_1(M))$ then we can show that
\[
[U\gamma U] = \coprod_{t \in T} U v_t,
\]
where~$T$ is a finite set, and~$v_t \in D^\times_f$. This means that we can write down the action of~$U_p$ or~$T_l$ on~$f$ as
\[
U_p f = [U\gamma U]f = \sum_{t \in \T} f|_k v_t,
\]
where~$v_t \in \A_f$ and~$u_{t,p}$ is the value of~$u_t$ at the place~$p$.

These~$v_t$ are explicitly computable; because we know that~$\gamma$ is trivial at all but finitely many places, we know that~$v_t$ must also be trivial at all but that finite set of places. In particular, let~$G=\left\{\abcds \in \GL_2(\Z_p): \abcds \equiv \smatrix{*}{*}{0}{1}\right\}$. If we are considering the automorphic forms used in Theorem~\ref{standard-jacquet-langlands}, then we can restrict to considering the double coset
\begin{equation}
\label{up-decomposition}
G\smatrix{p}{0}{0}{1}G = \coprod_{i=1}^p G \smatrix{p}{0}{pi}{1}
\end{equation}
for the operator~$U_p$, and
\begin{equation}
\label{tl-decomposition}
G\smatrix{\ell}{0}{0}{1}G = \coprod_{i=1}^p G \smatrix{1}{i}{0}{\ell} \coprod G\smatrix{\ell}{0}{0}{1}
\end{equation}
for the operator~$T_\ell$.

We can use these coset representatives to give us the~$v_t$; these are trivial, so equal to~1, at all places not equal to~2 or~$p$, and are given by the matrices in~\eqref{up-decomposition} or~\eqref{tl-decomposition}. Given these, we can then use the fact that we can decompose every element of~$D^\times_f$ into a unique product~$d \cdot d_i \cdot u$, where~$d \in D^\times$, $d_i$ is one of the coset representatives and~$u \in U$ (the open compact subgroup of~$D^\times_f$) and the definition of the right action in~\eqref{right-action} to rewrite the Hecke operator action as
\[
[U\gamma U]f = \sum_{t=1}^p f(d_i)|[u_{t,p}^{-1}]_k,
\]
where~$v_t=d \cdot d_i \cdot u_t$.

We note that for later sections, we will need to consider both the value at~$p$ and the value at~2; this will affect the decompositions given. This corresponds to choosing different subgroups of~$O^\times_D$; we will be using the intersection of~$O^\times_D$ with the subring of the quaternions~$D_{\Q_2}$ which have coefficients in~$\Z_2$ (the Lipschitz integers).

It can be shown that both the~$U_p$ operator acting on overconvergent modular forms and the~$U_p$ operator acting on automorphic forms are \emph{compact}; this means that it has a well-defined characteristic power series, which we can approximate by a series of finite-rank matrices, following~\cite{serre}. This means that we can effectively compute approximations to it using \Pari{}~\cite{Pari-tutorial}. 

We will also find it useful in practice to be able to consider smaller spaces of modular forms than the whole space of forms; this will make our computations run faster and enable us to consider higher weights and levels.
\subsection{Checking our computations}
If we are computing spaces of classical automorphic forms $p$-adically, then because these are finite-dimensional we can actually compute the action of Hecke operators on them exactly, because we know that their characteristic power series have integral coefficients of bounded size. We compute them to a sufficiently high $p$-adic precision, and then we can rewrite these $p$-adic numbers as rational integers.

Let us consider some numerical examples. Firstly, we compute the Hecke polynomial of~$T_3$ acting on the space~$S_5^D(U_1(7))$ of classical automorphic forms; it is
\[
(x^4 + 288x^2 + 20448)\cdot(x^4 + 18x^3 + 39x^2 - 1242x + 4761).
\]
It can be verified (using \Magma{}~\cite{magma}, say) that the Hecke 
polynomial 
of~$T_3$ acting on the space of classical modular forms~$S_5^{2-\rm 
new}(\Gamma_1(7)\cap\Gamma_0(2))$ is the same as this.

Again using \Magma{}, we can compute the Hecke polynomial of~$U_{11}$ 
acting on 
the spaces~$S_3^{2-\rm new}(\Gamma_1(11)\cap\Gamma_0(2))$ and~$S_3^{4- \rm new}(\Gamma_1(11))\cap\Gamma_0(4))$ of classical modular forms, and we see that these polynomials are~$x^2-14x+121$ and~$x^2+22x+121$, respectively. Using \Pari{}, we see that the Hecke polynomial of~$U_{11}$ on the space of automorphic forms~$S^D_3((1+\m)\cdot U_1(11))$ is given by
\[
(x^2-14x+121)\cdot(x^2+22x+121)^2,
\]
as predicted (we see the part which comes from 4-new forms appearing with multiplicity two here).

We can also compute automorphic forms with more restrictive level structure at~2. The characteristic polynomial of the Hecke operator~$U_5$ acting on the space of automorphic forms of weight~2 and level~$\Gamma_0(5)$ with level structure~$1+\m^3$ at~2 is given by
\[
(x^2-4x+5)^4\cdot(x^2+2x+5)^2\cdot(x^2-2x+5)^3\cdot(x-1)^3\cdot(x+1)^2\cdot(x-5),
\]
where the~$x-5$ comes from the norm form (because~$k=2$), the $(x-1)^3$ and the~$x^2+2x+5$ factors come from classical modular forms of level~$\Gamma_0(40)$, the~$(x+1)^2$ from modular forms of level~$\Gamma_0(20)$, and the quadratic factors~$(x^2-4x+5)^4$ and~$(x^2-2x+5)^3$ come from classical modular forms of level~$\Gamma_0(80)$.

One important feature of these calculations is that they give us an independent way to check that the modular forms algorithms in \Magma{} are giving the correct answers. Those algorithms use the theory of modular symbols and were programmed independently to the current work, so because we are using two different algebra packages and two different sets of programs, and still getting the same answer, we can be more certain that our programs are giving the correct results.
\section{Jacquet-Langlands in weight~1}

In this section we will prove the main theorem of our paper. Firstly, we introduce the concept of slopes of modular forms, which will be useful for us.

Let~$f$ be a normalized eigenform; either a classical modular form or an overconvergent modular form. We define the \emph{slope} of~$f$ to be the normalized $p$-valuation of the eigenvalue of~$U_p$ acting on~$f$.

\begin{theorem}[Coleman~\cite{coleman-overconvergent}, Theorem~1.1]
Let~$f$ be a classical modular eigenform of weight~$k$. Then the normalized $p$-slope of~$f$ is less than or equal to~$k-1$.

Conversely, if~$f$ is a $p$-adic overconvergent modular form of weight~$k$ with normalized slope strictly less than~$k-1$, then~$f$ is a classical modular form.
\end{theorem}
We see that there is a slight asymmetry in this result; if an overconvergent modular eigenform of weight~$k$ has normalized slope exactly~$k-1$, then it can be either classical or non-classical. There are examples of both; we will see this in weight~1 in Section~\ref{approximation}. The question of telling whether an overconvergent form of weight~$k$ and slope~$k-1$ is classical or not is raised by Coleman in~\cite{coleman-old-jntb}; see Section~7, and is still open in general.

We now prove the main theorem of this paper; that the standard Jacquet-Langlands correspondence can be extended to weight~1.
\begin{theorem}
\label{main-theorem}
Let~$N$ be a positive odd integer and let~$i$ be either~0 or~1.


If~$f \in S_1^{2^i-\text{new}}(\Gamma_1(N) \cap \Gamma_0(2^{i+1}))$, then there exists an overconvergent automorphic form~$f_A \in S^{D,\dagger}_1(U_1(N)\cdot G)$ with the same Hecke eigenvalues as~$f$ (if~$i=0$, then there is no extra level structure at~2, and if~$i=1$ then~$G=1+\m$). Conversely, if~$f_A$ is an overconvergent automorphic form of weight~1, then there exists an overconvergent modular form~$f$ of weight~1 with the same Hecke eigenvalues as~$f_A$.

\end{theorem}
We note that a version of this is true in more generality, for other subgroups of~$D^\star$ of finite index, but we will not need this for the section on approximation eigenforms.
\begin{proof}
Using the theory of families of modular forms (for instance, \cite{wan} proves the existence of suitable families of modular forms), we can find a $p$-adic family~$\{f_i\}_{i \in \N}$ of classical modular forms with each~$f_i$ having weight~$1+(p-1)\cdot p^i$, with each~$f_i$ having Fourier expansion congruent to that of~$f$ modulo~$p^{i+1}$. We can use the Jacquet-Langlands correspondence to find a family of classical automorphic forms~$\left\{f_{A_i}\right\}_{i \in \N}$, each of which has the same Hecke eigenvalues as~$f_i$, and then we see that if we take the limit of these~$f_{A,i}$ then it is the overconvergent automorphic form~$f_A$ of weight~1.

Conversely, if we take an overconvergent automorphic form~$f_A$ of weight~1 and slope~0, then we can fit it into a $p$-adic family of classical automorphic forms~$\{f_{A,i}\}_{i \in \N}$ of weight~$1+(p-1)\cdot p^i$, each of which has Hecke eigenvalues which are congruent to those of~$f_A$ modulo~$p^{i+1}$. We then use the Jacquet-Langlands correspondence in the other direction to find a family of classical modular forms~$\left\{f_i\right\}_{i \in \N}$, and then we take the limit of these to find~$f$.

We note that in the recent work of Chenevier~\cite{chenevier}, we see that if we have a $p$-adic family of automorphic forms, then we can apply the classical Jacquet-Langlands correspondence to show that the image of this is a $p$-adic family of modular forms.
\end{proof}
We see that to prove this result we had to use the full force of the Jacquet-Langlands correspondence in both directions. We can partially relax this by using results from Section~2.5 of~\cite{hida-hilbert-book}, which use the theory of theta series to give a map from automorphic forms to modular forms without using the Jacquet-Langlands correspondence.

This works in the following way: given a classical automorphic form~$f_A$, we can create a theta series~$\theta_{f_A}$ which is a classical modular form and has the same Hecke eigenvalues as~$f_A$. We then create a family of classical automorphic forms which approximate the overconvergent automorphic form of weight~1, and then the family of theta series will approximate an overconvergent modular form of weight~1 (which may or may not be classical).
\section{Approximating eigenforms}
\label{approximation}

In this section we will give an account of how to actually find approximations to overconvergent automorphic eigenforms of weight~1, using \Pari{}  programs. We also indicate how this method can be generalized to find other forms.

This method is a development of the work of \Gouvea{} and Mazur in~\cite{gouvea-mazur-searching}, where they find overconvergent 5-adic modular eigenforms of weight~0 by iterating the action of the~$U_5$ operator. This in turn builds on the work of Atkin and O'Brien~\cite{atkin-obrien} which pioneered this technique for finding $p$-adic eigenforms for~$p=13$. 

On the modular side, we will consider the space of classical modular forms $S^{4-\rm new}_1(\Gamma_1(11)\cap\Gamma_0(4))$; this can be checked to be one-dimensional, and it is in fact generated by the $\eta$-product $f:=\eta(q^{2})\eta(q^{22})$, which is necessarily a Hecke eigenform. This has Fourier expansion at~$\infty$ given by
\[
f(q)=q\prod_{n=1}^\infty (1-q^{2n})(1-q^{22n})=q - q^3 - q^5 + q^{11} + q^{15} + O(q^{23});
\]
in particular, it has 11-slope~0, which shows that it is in the interesting case left open by the theory of Coleman, where the slope is~$k-1$. By the Ramanujan-Petersson Conjecture, the Fourier coefficients~$a_p$ of~$f(q)$ satisfy~$|a_p| < 2$.

We use code based on the work of Jacobs outlined in Section~\ref{computations} to find the slopes of the~$U_{11}$ operator acting on weight~1 overconvergent automorphic forms of level~$U_0(11)\cdot(1+\mathfrak{m})$. This tells us that the lowest slopes are~$0,0,1,2,2,2$ with that multiplicity.

Because we know from Theorem~\ref{jacquet-langlands-extended-classical} and Theorem~\ref{main-theorem} that there should be two automorphic forms related to the classical modular form~$f$ appearing on the automorphic side, and from above we see that~$f$ has slope~0, we see that the two slope~0 automorphic forms are the ones that correspond to~$f$, and all of the others do not correspond to classical forms.

Let~$\{h_i\}$ be the set of simultaneous eigenforms for~$S^D_1(U_1(11)\cdot(1+\mathfrak{m}))$, so let~$g_1$ be a random nonzero element of this space. We can think of~$g_1$ as a triple of power series in~$z$, because we know that an automorphic form is determined by its values on the tuple~$\{d_i\}$. (In fact, to aid the calculations we can just choose a triple of polynomials, and make two of them zero).

We assume that we can write~$g_1$ as a linear combination of these eigenforms:
\[
g_1=\alpha  h_0 + \beta h_1 + \gamma h_2 + \cdots,
\]
where~$h_0$ and~$h_1$ are the two slope~0 forms.
We now compute the action of~$U_{11}^N$ on~$g_1$, for some large integer~$N$. We see that $U_{11}^N(g_1)$ will be congruent to~$\alpha h_0 + \beta h_1$ modulo~$p^N$, because the action of~$U_{11}^N$ on the other~$h_i$ will include a multiplication by~$p^{N}$ at least. This means that we have an approximant to the sum of the eigenvectors~$h_0$ and~$h_1$.

We now choose a second random element~$g_2$ and compute~$U_{11}^N(g_2)$. This will also be congruent to a linear combination of~$h_0$ and~$h_1$ modulo~$p^N$; with very high probability, these two linear combinations are linearly independent, and we can now use linear algebra to find~$h_0$ and~$h_1$ modulo~$p^N$ from them.

Finally, to find eigenforms for \emph{all} of the Hecke operators, we consider the action of the~$W$ operator, the analogue of a diamond operator in the classical setting, which is defined to be
\[
Wf=[U(1+i)U]f.
\]
The action of this operator splits the 2-dimensional eigenspace for~$U_{11}$ into two one-dimensional eigenspaces. Basis elements for each of these eigenspaces are eigenforms for all of the Hecke operators~$T_\ell$ (for $\ell$ a prime not equal to~2 or~11) and~$U_{11}$.

It would be interesting if one could find simultaneous eigenforms for~$U_p$ and for the other Hecke operators exactly (rather than approximately) by a similar process, given the eigenvalues. In this example, we know that such an eigenform has~$U_{11}$-eigenvalue~1, so it satisfies an equation of the form
\[
f(z)=1\cdot f(z) = (U_{11}f)(z)=\sum_{i=0}^{11} f(\gamma_i z) \cdot (c_iz+d_i)^{-1},
\]
where the~$\gamma_i=\abcdsi$ are a set of matrices that represent the Hecke operator~$U_{11}$. We can write down similar recurrence relations for each of the Hecke operators. It would be useful to be able to solve these explicitly.

This process can be performed in more generality, to find approximations to eigenforms of slope~0 and of higher slope. The higher-slope cases are more delicate; we cannot use the same techniques as above in general because we need to divide by powers of~$p$, which will reduce the accuracy at which we are working. The same methods do work for slope~1 in the example we have considered here, because we gain more accuracy at each step by iterating the~$U_{11}$ operator than we lose by dividing by~$11$.

For higher slopes, a generalization of the methods used by Loeffler~\cite{loeffler}, Section~5, would seem more appropriate. These use the properties of an inner product on the space of 5-adic overconvergent modular forms to compute spectral expansions of eigenfunctions for the~$U_5$ operator.


We note here that the methods we have outlined will also work for higher weight forms; let~$k$ be a positive integer. We can find any automorphic forms of slope~0 using exactly this procedure; these will be classical automorphic forms, so they will be determined by a tuple of polynomials. After subtracting these out, we will be able to find forms of higher slope, and this will enable us to approximate overconvergent automorphic forms of weight~$k$.

\section{Telling classical from non-classical automorphic forms}

It is well-known that $p$-adic modular forms of weight~$k$ and slope exactly~$k-1$ can be either classical or non-classical. Examples of both can be exhibited; the form~$f$ from the previous section was classical, and we exhibit an example of a non-classical form below.

Using the same approximation techniques as in Section~\ref{approximation}, we can find an approximation to an overconvergent automorphic form of slope~0, weight~1 and level~$U_1(7)$ (with no level condition at~2). We know that such a form must exist because we can compute the slopes of~$U_7$ acting on weight~1 overconvergent automorphic forms, and we find that these are~$0,1,1,2,2,2$. This means that there is a unique form~$g$ (which aids the calculations) with eigenvalue~$1 + 5\cdot7 + 4\cdot7^2 + 5\cdot7^3+O(7^4) $. This is \emph{not} a classical modular form, because we can check that there are no classical modular cuspforms of weight~1 and level~$\Gamma_1(7)\cap\Gamma_0(2)$.

Using work on families of modular forms we can find classical automorphic forms of weight~$1+(p-1)p^n$, for any non-negative integer~$n$, which lie in the same Hida family as our weight~1 form, whether or not it is classical; in fact, because the slope is~0, using the theorem of Coleman quoted above we can see that the only form in the family for which there is a doubt as to whether it is classical is the weight~1 form.

It would be interesting to have a method for telling an overconvergent automorphic form which comes from a classical modular form (via a generalized Jacquet-Langlands correspondence) from an overconvergent automorphic form which does \emph{not} come from a classical modular form. 

In some circumstances, we may know, following the version of the classical Jacquet-Langlands correspondence given as Theorem~\ref{jacquet-langlands-extended-classical}, that a classical modular form will appear with multiplicity two on the automorphic side, but we would like to have a more intrinsic criterion. Also, we may not be given all of the automorphic forms; we would like a method that only requires us to consider one form.

If we are willing to invoke both a Jacquet-Langlands correspondence \emph{and} the Ramanujan-Petersson conjecture on coefficients of classical modular forms, then we can show that certain automorphic forms are \emph{not} coming from a classical modular form. The Ramanujan-Petersson conjecture tells us that, if~$f(q)=\sum_{n \in \N} a_n q^n$ is the Fourier expansion of a normalized cuspidal modular eigenform of weight~$k$, and~$\ell$ is a prime, then
\[
|a_\ell | \le 2 p^{\frac{k-1}2}.
\]
This was conjectured by Ramanujan, generalized by Petersson, and proved by Deligne~\cite{deligne-weil-conjectures} as a consequence of the Weil conjectures.

We can approximate the eigenvalues of the Hecke operators acting on automorphic forms to a high degree of accuracy using our computer programs, and if we can show that these eigenvalues are not algebraic numbers, then we are done. We can bound the maximum degree that these algebraic numbers can have by dimension considerations. We can use this criterion to show that the weight~1 form at level~$U_1(7)$ is not classical; however, we are using several very high-level results to obtain this.

\section{Acknowledgements}
This paper was written while I was the GCHQ Research Fellow at Merton College, Oxford and a Research Fellow at the University of Bristol. I would like to thank the College and the University for their hospitality during my terms.

I would like to thank Kevin Buzzard for many helpful conversations and for much help and guidance with this project. 
Without his guidance this paper could not have been written. I would also like to thank Gabor Wiese for helpful 
conversations. All remaining errors are of course due to the author.

I would also like to thank Daniel Jacobs for his help and for his work on classical and overconvergent automorphic forms.

I would also like to thank William Stein for giving me an account on his machine {\sc Meccah}, on which many of the computations were performed.


\end{document}